\documentclass[12pt]{imsart}
\usepackage{amsthm,amsmath}
\usepackage[colorlinks,citecolor=blue,urlcolor=blue]{hyperref}

\usepackage{amsmath,amsfonts,amssymb,amsthm}

\usepackage{enumitem}

\RequirePackage{amsthm,amsmath,amsfonts,amssymb}
\RequirePackage[numbers]{natbib}
\RequirePackage[colorlinks,citecolor=blue,urlcolor=blue]{hyperref}
\RequirePackage{graphicx}
\usepackage{lineno,hyperref}

\startlocaldefs
\modulolinenumbers[5]

\usepackage[utf8]{inputenc}
\usepackage[english]{babel}

\usepackage[margin=1in]{geometry}
\usepackage{bm}
\usepackage{graphicx}
\usepackage{amssymb,amsmath,amsthm,mathrsfs}
\usepackage{epstopdf}
\usepackage{algorithm}
\usepackage{amsfonts,delarray}
\usepackage{algorithm,algcompatible,amsmath}
\usepackage{rotating,color}
\usepackage{subfigure}
\usepackage{multirow}
\usepackage{titlesec,textcase}
\usepackage{longtable}
\usepackage{bbm}
\usepackage{rotating}

\newtheorem{Theorem}{Theorem}[section]

\newtheorem{Corollary}[Theorem]{Corollary}

\theoremstyle{definition}

\RequirePackage[normalem]{ulem} 

\definecolor{rp}{RGB}{83,54,106}

\def\boxit#1{\vbox{\hrule\hbox{\vrule\kern6pt\vbox{\kern6pt#1\kern6pt}\kern6pt\vrule}\hrule}}

\endlocaldefs

\begin{document}
\begin{frontmatter}
\title{On the Randi\'{c} index and its variants of network data}

\runtitle{Randi\'{c} index of network data}
\runauthor{ }
\begin{aug}

\author[A]{\fnms{Mingao} \snm{Yuan}\ead[label=e1]{mingao.yuan@ndsu.edu}}



\address[A]{Department of Statistics,
North Dakota State University,
\printead{e1}}
\end{aug}

\begin{abstract}
 Summary
statistics play an important role in network data analysis. They can provide us with meaningful insight into the structure of a network. The Randi\'{c} index is one of the most popular network statistics that has been widely used for quantifying information of biological networks, chemical networks,  
pharmacologic networks, etc. A topic of current interest is to find bounds or limits of the Randi\'{c} index and its variants. A number of bounds of the indices are available in literature. Recently, there are several attempts to study the limits of the indices in the Erd\H{o}s-R\'{e}nyi random graph by simulation.
In this paper, we shall derive the limits of the Randi\'{c} index and its variants of an inhomogeneous Erd\H{o}s-R\'{e}nyi random graph. Our results charaterize how network heterogeneity affects the indices and provide new insights about the Randi\'{c} index and its variants.  Finally we apply the indices to several real-world networks.
\end{abstract}



\end{frontmatter}

\section{Introduction}
\label{S:1}

 A network (graph) consists of a set of agents and a set of pairwise interactions
among the agents. Networks are canonical models that capture relations within or between data sets. Due to the increasing popularity of relational data, network data analysis has been a primary research topic in statistics, machine learning and many other scientific fields \cite{BS16,A18,K09,N09,GZFA10}.
One of the fundamental problems in network data analysis is to understand the structural properties of a given network. The structure of a small network can be easily described by its visualization. However, larger networks can be difficult to envision and describe. It is thus important to have  several summary statistics that provide us with meaningful insight into the structure of a network. Based on these statistics, we are able to compare networks or classify
them according to properties that they exhibit.
There are a wealth of descriptive statistics that measure some aspect of the structure or characteristics of a network. For example, the diameter of a network measures the maximum distance between two individuals; the global clustering coefficient measures the extent to which individuals in a graph tend to cluster together; the modularity is a measure of the strength of division of a network into subgroups.

Summary statistics of  networks are sometimes termed topological indices, especially in chemical or pharmacological science \cite{MCSGDF18}. 
One of the most popular topological indices is the Randi\'{c} index invented in 
\cite{R75}. The Randi\'{c} index measures the extent of branching of a network \cite{BT78,R75}.  It was observed that the Randi\'{c} index is strongly correlated with a variety of physico-chemical properties of alkanes \cite{R75}. The Randi\'{c} index play a central role in understanding quantitative structure-property and structure-activity relations in chemistry and pharmocology \cite{R08,RNP16}. In subsequent years,  the Randi\'{c} index finds countless applications. For instance, it is used to characterize and quantify the similarity between different networks or subgraphs of the same network \cite{FT13}, it serves as a quantitative characterization of network heterogeneity \cite{E10}, and graph robustness can be easily estimated by the Randi\'{c} index \cite{DMRSV18,DMMRSF21}. Moreover, the Randi\'{c} index possesses a wealth of non-trivial and interesting mathematical properties \cite{BE98,BES99,CFK10,DSG17,LS08}. Motivated by the Randi\'{c} index, various Randi\'{c}-type indices have been introduced and attracted great interest in the past years. Among them, the harmonic index is a well-known one \cite{F87,FMS93, Z12,RS17}.

One of the popular research topics in the study of topological indices is to derive bounds of
the indices and study their asymptotic properties. Recently, \cite{MMRS20,MMRS21} performed numeric and analytic analyses of the Randi\'{c} index and the harmonic index in  the Erd\H{o}s-R\'{e}nyi random graph. Analytic upper and lower bounds of the two indices are obtained and simulation studies show that the indices converge to one half of the number of nodes. Additionally, \cite{DMRSV18,DHHIR20,LSG21} find the expectations of  variants of the Randi\'{c} index in the Erd\H{o}s-R\'{e}nyi random graph. 
However, these results only apply to the Erd\H{o}s-R\'{e}nyi random graph and the exact limits of the indices are not theoretically studied.

In this paper, we shall derive the limits of the general Randi\'{c} index and the general sum-connectivity index in an inhomogeneous Erd\H{o}s-R\'{e}nyi random graph. The general Randi\'{c} index and the general sum-connectivity index contain the Randi\'{c} index and the harmonic index as a special case, respectively. Thus our results theoretically validate the empirical observations in \cite{MMRS20,MMRS21} that the indices of the Erd\H{o}s-R\'{e}nyi random graph converge to one half of the number of nodes. In addition, our results explicitly describe how network heterogeneity affects the indices. We also observe that the limits of the Randi\'{c} index and the harmonic index do not depend on the sparsity of a network, while the limits of their variants do. In this sense, the Randi\'{c} index and the harmonic index are more preferable than their variants as measures of network structure.

The structure of the article is as follows. In Section \ref{main} we present the main results. Section \ref{simu} summarizes simulation results and real data application. The proof is deferred to Section \ref{proof}.

\medskip

Notations:  Let $c_1,c_2$ be positive constants and  $n_0$ be a positive integer. For two positive sequence $a_n$, $b_n$, denote $a_n\asymp b_n$ if $c_1\leq \frac{a_n}{b_n}\leq c_2$  for $n\geq n_0$; denote $a_n=O(b_n)$ if $\frac{a_n}{b_n}\leq c_2$ for $n\geq n_0$; $a_n=o(b_n)$ if $\lim_{n\rightarrow\infty}\frac{a_n}{b_n}=0$. Let $X_n$ be a sequence of random variables. $X_n=O_P(a_n)$ means $\frac{X_n}{a_n}$ is bounded in probability.  $X_n=o_P(a_n)$ means $\frac{X_n}{a_n}$ converges to zero in probability. Denote $a_+=\max\{a,0\}$.

\section{The Randi\'{c} index and its variants}\label{main}

A graph is a mathematical model of network that consists of nodes (vertices) and edges. Let $\mathcal{V}=[n]:=\{1,2,\dots,n\}$ for a given positive integer $n$.
An \textit{undirected} graph on $\mathcal{V}$ is a pair $\mathcal{G}=(\mathcal{V},\mathcal{E})$ in which $\mathcal{E}$ is a collection of subsets of $\mathcal{V}$ such that $|e|=2$ for every $e\in\mathcal{E}$. Elements in $\mathcal{E}$ are called edges. A graph can be conveniently represented as an adjacency matrix $A$, where $A_{ij}=1$ if $\{i,j\}$ is an edge, $A_{ij}=0$ otherwise and $A_{ii}=0$. It is clear that $A$ is symmetric, since $\mathcal{G}$ is undirected. A graph is said to be random if $A_{ij} (1\leq i<j\leq n)$ are random.

 Let $f=(f_{ij})$, $(1\leq i<j\leq n)$ be a vector of numbers between 0 and 1.
The inhomogeneous Erd\H{o}s-R\'{e}nyi random graph $\mathcal{G}(n,p_n, f)$ is defined as
\[\mathbb{P}(A_{ij}=1)=p_n f_{ij},\]
where $p_n\in[0,1]$ and $A_{ij}\ (1\leq i<j\leq n)$ are independent. If all $f_{ij}$ are the same, then $\mathcal{G}(n,p_n, f)$ is the Erd\H{o}s-R\'{e}nyi random graph. For a non-constant vector $f$, $\mathcal{G}(n,p_n, f)$ is an inhomogeneous version of the Erd\H{o}s-R\'{e}nyi random graph. This model covers several random graphs that have been extensively studied in random graph theory and algorithm analysis \cite{CHHS20,CCH20,CHHS21,CGL16,YXL21}.

Given a constant $\alpha$, the general Randi\'{c} index of a graph $\mathcal{G}$ is defined as (\cite{BE98}) 
\begin{equation}\label{randicindex}
\mathcal{R}_{\alpha}=\sum_{\{i,j\}\in \mathcal{E}}d_i^{\alpha}d_j^{\alpha}, 
\end{equation}
where $d_k$ is the degree of node $k$, that is, $d_k=\sum_{j\neq k}A_{kj}$. The index $\mathcal{R}_{\alpha}$ generalizes the well-known Randi\'{c} index $\mathcal{R}_{-\frac{1}{2}}$ invented in \cite{R75}. When $\alpha=-1$, the index $\mathcal{R}_{-1}$ corresponds to the modified second Zagreb index \cite{NKMT03,CFK10}. 

Another popular variant of the Randi\'{c} index is the general sum-connectivity index \cite{ZT09,ZT10} defined as
\begin{equation}\label{harindexm}
\chi_{\alpha}=\sum_{\{i,j\}\in \mathcal{E}}(d_i+d_j)^{\alpha}.
\end{equation}
An important special case is the harmonic index $\mathcal{H}=2\chi_{-1}$ \cite{F87,FMS93, Z12}.

Recently,
\cite{MMRS20,MMRS21} conduct a simulation study of the Randi\'{c} index $\mathcal{R}_{-\frac{1}{2}}$ and the harmonic index $\mathcal{H}=2\chi_{-1}$ in the Erd\H{o}s-R\'{e}nyi random graph and observe that the indices converge to $n/2$. Moreover,  \cite{DMRSV18,DHHIR20,LSG21} derive analytical
expressions of the expectations for the indices $\mathcal{R}_{-1}$,$\mathcal{\chi}_{1}$,$\mathcal{\chi}_{2}$ of the Erd\H{o}s-R\'{e}nyi random graph. In this paper,
we shall derive the exact limits of the general Randi\'{c} index $ \mathcal{R}_{\alpha}$ and the general sum-connectivity index $\chi_{\alpha}$ in $\mathcal{G}(n,p_n,f)$. Our results significantly improve the results in \cite{DMRSV18,MMRS20,MMRS21,DHHIR20,LSG21} and provide new insights about the Randi\'{c} index and its variants.

\begin{Theorem}\label{theorem:1} Let $\alpha$ be a fixed constant and $\mathcal{G}(n,p_n,f)$ be the inhomogeneous Erd\H{o}s-R\'{e}nyi random graph. Suppose $ np_n\log 2\geq \log n$ and $\min_{1\leq i<j\leq n}\{f_{ij}\}>\epsilon$ for some positive constant $\epsilon\in(0,1)$. Then
\begin{eqnarray}\label{randic0}
   \mathcal{R}_{\alpha}&=&\left[1+O_P\left(\frac{(\log(np_n))^{4(1-\alpha)_+}}{\sqrt{np_n}}\right)\right]
  p_n^{2\alpha+1}\sum_{i<j}f_i^{\alpha}f_j^{\alpha}f_{ij},\\
  \label{harmonic0}
   \mathcal{\chi}_{\alpha}&=&\left[1+O_P\left(\frac{(\log(np_n))^{2(1-\alpha)_+}}{\sqrt{np_n}}\right)\right]p_n^{\alpha+1}\sum_{i<j}(f_i+f_j)^{\alpha}f_{ij},
\end{eqnarray}
where $f_i=\sum_{j\neq i}^nf_{ij}$. 
\end{Theorem}
\medskip
 
The condition $\min_{1\leq i<j\leq n}\{f_{ij}\}>\epsilon$ implies the minimum expected degree scales with $np_n$. The condition $ np_n\log 2\geq \log n$ means that the graph is relatively dense. A similar condition is assumed in \cite{CHHS20} to study the maximum eigenvalue of the inhomogeneous random graph.

Note that the expected total degree of $\mathcal{G}(n,p_n,f)$ has order $n^2p_n$. Thus $p_n$ controls the sparsity of the network: a graph with smaller $p_n$ would have fewer edges. By (\ref{randic0}) and (\ref{harmonic0}), the limits of the Randi\'{c} index $\mathcal{R}_{-\frac{1}{2}}$ and the harmonic $\mathcal{\chi}_{-1}$ do not depend on $p_n$, while the limits of their variants do involve $p_n$. Asymptotically, the Randi\'{c} index and the harmonic are uniquely determined by the network structure parametrized by $f$. 
In this sense, they are superior to their variants as measures of global structure of networks.

Now we present two examples of $\mathcal{G}(n,p_n,f)$.
The simplest example is the Erd\H{o}s-R\'{e}nyi random graph, that is, $f_{ij}\equiv 1$. We denote the graph as $\mathcal{G}(n,p_n)$. 

\begin{Corollary}\label{cor:1} Let $\alpha$ be a fixed constant. For the Erd\H{o}s-R\'{e}nyi random graph $\mathcal{G}(n,p_n)$ with $ np_n\log 2\geq \log n$, we have
\begin{eqnarray}\label{randicer}
   \mathcal{R}_{\alpha}&=&
   \frac{n^{2(1+\alpha)}p_n^{2\alpha+1}}{2}\left[1+O_P\left(\frac{(\log(np_n))^{4(1-\alpha)_+}}{\sqrt{np_n}}\right)\right],\\ \label{randicer2}
   \mathcal{\chi}_{\alpha}&=&2^{\alpha-1}n^{\alpha+2}p_n^{\alpha+1}\left[1+O_P\left(\frac{(\log(np_n))^{2(1-\alpha)_+}}{\sqrt{np_n}}\right)\right].  
\end{eqnarray}
Especially, the Randi\'{c} index $\mathcal{R}_{-\frac{1}{2}}$ is equal to 
\[\mathcal{R}_{-\frac{1}{2}}=\frac{n}{2}\left[1+O_P\left(\frac{(\log(np_n))^{4(1-\alpha)_+}}{\sqrt{np_n}}\right)\right],\]
the modified second Zagreb  index $\mathcal{R}_{-1}$ is equal to
\[\mathcal{R}_{-1}=\frac{1}{2p_n}\left[1+O_P\left(\frac{(\log(np_n))^{4(1-\alpha)_+}}{\sqrt{np_n}}\right)\right],\] 
and the harmonic index $\mathcal{H}$ is equal to
\[\mathcal{H}=\frac{n}{2}\left[1+O_P\left(\frac{(\log(np_n))^{2(1-\alpha)_+}}{\sqrt{np_n}}\right)\right].\]
\end{Corollary}
\medskip

According to Corollary \ref{cor:1}, the ratio $ \frac{2}{n}\mathcal{R}_{-\frac{1}{2}}$ or $ \frac{2}{n}\mathcal{H}$ converges in probability to 1 when $ np_n\log 2\geq \log n$. This theoretically confirms the empirical observation in \cite{MMRS20,MMRS21} that the Randi\'{c} index  $\mathcal{R}_{-\frac{1}{2}}$ or the harmonic index $\mathcal{H}$ is approximately equal to  $\frac{n}{2}$.  
The expectation of the  indices  $\mathcal{R}_{-1},\mathcal{\chi}_{1},\mathcal{\chi}_{2}$ are derived in \cite{DMRSV18,DHHIR20,LSG21}. Our results show the indices are asymptotically equal to their expectations. Moreover, Corollary \ref{cor:1} clearly quantifies how $p_n$ affects the convergence rates: the larger $p_n$ is, the faster the convergence rates are.

In addition, (\ref{randicer}) and (\ref{randicer2}) explicitly characterize how the leading terms of $\mathcal{R}_{\alpha}$ and $  \mathcal{\chi}_{\alpha}$ depend on $\alpha$. Note that
\begin{eqnarray*}\label{errand}
\frac{n^{2(1+\alpha)}p_n^{2\alpha+1}}{2}&=&\frac{n}{2}(np_n)^{2\alpha+1},\\ \label{errhar}
2^{\alpha-1}n^{\alpha+2}p_n^{\alpha+1}&=&2^{\alpha-1}n(np_n)^{\alpha+1}.
\end{eqnarray*}
For given $n,p_n$ such that $ np_n\log 2\geq \log n$, the  leading terms are increasing functions of $\alpha$. The indices would be extremely large or small for large $|\alpha|$ and large $n$. In this sense, it is preferable to use $\mathcal{R}_{\alpha}$ or $  \mathcal{\chi}_{\alpha}$ with small $|\alpha|$ (for instance, $|\alpha|\leq1$).

Next, we provide a non-trivial example.
Let $f_{ij}=e^{-\kappa \frac{i}{n}}e^{-\kappa \frac{j}{n}}$ with a positive constant $\kappa$. Then $e^{-2\kappa}\leq f_{ij}\leq 1$ for  $0\leq i<j\leq n$. In this case, $\min_{1\leq i<j\leq n}\{f_{ij}\}>\epsilon$ holds with $\epsilon=e^{-2\kappa}$. Straightforward calculation yields $ f_i=ne^{-\kappa\frac{i}{n}}\frac{(1-e^{-\kappa})}{\kappa}(1+o(1))$ and
 \begin{eqnarray*}
 \sum_{i<j}f_i^{-1}f_j^{-1}f_{ij}&=&\frac{\kappa^2}{2(1-e^{-\kappa})^{2}}+o(1),\\
 \sum_{i<j}f_i^{\alpha}f_j^{\alpha}f_{ij}&=&\frac{n^{2(\alpha+1)}(1-e^{-\kappa})^{2\alpha}(1-e^{-(1+\alpha)\kappa})^2}{2(1+\alpha)^2\kappa^{2(\alpha+1)}}(1+o(1)),\ \ \alpha\neq-1,\\
 \sum_{i<j}(f_i+f_j)^{\alpha}f_{ij}&=&\frac{n^{\alpha+2}}{2}\left(\frac{1-e^{-\kappa}}{\kappa}\right)^{\alpha}\int_0^1\int_0^1\frac{\left(e^{-\kappa x}+e^{-\kappa y}\right)^{\alpha}}{e^{\kappa (x+y)}}dxdy+o(1).
 \end{eqnarray*}
 Then
\begin{eqnarray}\label{ran1}
   \mathcal{R}_{-1}&=&\left[1+O_P\left(\frac{(\log(np_n))^{2}}{\sqrt{np_n}}\right)\right] \frac{1}{2p_n} \frac{\kappa^2}{(1-e^{-\kappa})^{2}}
,\\ \label{ran2}
   \mathcal{R}_{\alpha}&=&\left[1+O_P\left(\frac{(\log(np_n))^{2}}{\sqrt{np_n}}\right)\right] \frac{n^{2(\alpha+1)} p_n^{2\alpha+1}}{2} \frac{ (1-e^{-\kappa})^{2\alpha}(1-e^{-(1+\alpha)\kappa})^2}{(1+\alpha)^2\kappa^{2(\alpha+1)}}
,\ \ \alpha\neq-1,\\
\label{ran3}
    \mathcal{\chi}_{\alpha}&=&\left[1+O_P\left(\frac{(\log(np_n))^{2}}{\sqrt{np_n}}\right)\right]\frac{n^{\alpha+2}p_n^{\alpha+1}}{2}\left(\frac{1-e^{-\kappa}}{\kappa}\right)^{\alpha}\int_0^1\int_0^1\frac{\left(e^{-\kappa x}+e^{-\kappa y}\right)^{\alpha}}{e^{\kappa (x+y)}}dxdy.
\end{eqnarray}

Since larger $\kappa$ makes the expected degrees more heterogeneous, the parameter $\kappa$ can be considered as heterogeneity level of the graph. As $\kappa$ increases,  $\mathcal{R}_{\alpha}$ or $\mathcal{\chi}_{\alpha}$ decreases if $\alpha>-1$, and $\mathcal{R}_{\alpha}$ or $\mathcal{\chi}_{\alpha}$ increases if $\alpha\leq-1$. This shows the effect of heterogeneity on $\mathcal{R}_{\alpha}$ or $\mathcal{\chi}_{\alpha}$. The indices could be used as indicators whether a network follows the  Erd\H{o}s-R\'{e}nyi random graph model.

\section{Real data application}\label{simu}

In this section, we apply the general Randi\'{c} index and the general sum index to the following real-world networks: `karate', `macaque', `UKfaculty', `enron', `USairports', `immuno', `yeast'. These networks are available in  the `igraphdata' package of R. 

For each network, the indices $\mathcal{R}_{-\frac{1}{2}}$, $\mathcal{R}_{-1}$, $\mathcal{\chi}_{-\frac{1}{2}}$, $\mathcal{\chi}_{-1}$ and the bound $\log n/(n\log 2)$ are calculated.  Here, $\log n/(n\log 2)$ is the sparsity lower bound required by Theorem \ref{theorem:1} and Corollary \ref{cor:1}. In addition, we also compute several descriptive statistics: the number of nodes ($n$), the edge density, the maximum degree ($d_{max}$), the median degree ($d_{mean}$) and the minimum degree ($d_{min}$). 
These results are summerized in Table \ref{tabreal}. The edge densities of networks `macaque', `UKfaculty', `enron' and `USairports'  are greater than  $\log n/(n\log 2)$, which indicates our theoretical results are applicable.  The Randi\'{c} indices $\mathcal{R}_{-\frac{1}{2}}$ and the harmonic indices $2\mathcal{\chi}_{-1}$ of `enron' and `USairports'  are much smaller than $\frac{n}{2}$, the indices of the Erd\H{o}s-R\'{e}nyi random graph. Thus the Erd\H{o}s-R\'{e}nyi random graph may not be a good model for these two networks.  The networks `macaque' and `UKfaculty' have the indices close to $\frac{n}{2}$. In this sense, they can be considered as samples from the Erd\H{o}s-R\'{e}nyi random graph model. For the networks `karate', `immuno' and `yeast',  the edge densities are slightly smaller than the bound $\log n/(n\log 2)$. Note that the condition $p_n>\log n/(n\log 2)$ is a sufficient condition for Theorem \ref{theorem:1} and Corollary \ref{cor:1} to hold and can not be relaxed  based on the current proof technique. We conjecture that Theorem \ref{theorem:1} and Corollary \ref{cor:1} still hold if $np_n\rightarrow\infty$. Currently, we are not clear whether our theoretical results can be applied to the networks `karate', `immuno' and `yeast' or not.For sparse networks, that is, $np_n=O\left(1\right)$, the Randi\'{c} index $\mathcal{R}_{-\frac{1}{2}}$ could assume any value between 0 and $\frac{n}{2}$, which is empirically verified in \cite{MMRS20}. Therefore, the Randi\'{c} index $\mathcal{R}_{-\frac{1}{2}}$ far  less than $\frac{n}{2}$ does not necessarily imply the network are not generated from the Erd\H{o}s-R\'{e}nyi random graph model.   We point out that a statistical hypothesis testing is needed to test whether the Randi\'{c} index is equal to some number. Based on our knowledge, there is no such test available in literature. It is an interesting future topic to propose a test for the Randi\'{c} index.

\begin{table}[h]
\begin{center}
\begin{tabular}{ |c|c|c| c|c|c|c|c|c|c|c|} 
 \hline
 network & $n$& $\log n/(n\log 2)$ & density & $d_{max}$&$d_{median}$& $d_{min}$ & $\mathcal{R}_{-\frac{1}{2}}$ & $\mathcal{R}_{-1}$ & $\mathcal{\chi}_{-\frac{1}{2}}$ & $\mathcal{\chi}_{-1}$ \\  
  \hline
   karate     & 34 &0.149 & 0.134  & 17 & 5 &  3 &13.970&2.866&  21.001 & 5.927\\ 
  \hline

    macaque     & 45 & 0.122&  0.251  & 22  & 11  & 4 & 21.576 & 2.092 & 50.702& 10.374 \\
      \hline
     UKfaculty     & 81 & 0.078 & 0.175 & 41  & 13 & 2  & 37.728 & 2.957 & 99.101 & 17.738 \\ 
  \hline
     enron     & 184 &0.040 & 0.130  & 111  & 31 &21  &80.876 &4.063& 276.792 & 37.672\\ 
  \hline
    USairports     & 755&0.012 & 0.016  & 168& 11 &5  &262.836&41.776& 602.894    &106.592\\ 
  \hline
       immuno     & 1316 & 0.0078  & 0.0072   &  17  & 10  & 3   & 648.820 & 70.951 & 1410.842 & 320.022 \\ 
  \hline
     yeast     & 2617& 0.004 & 0.003  & 118 &10 & 4  &1076.274 & 285.491&  2034.479 &469.020\\ 
  \hline
\end{tabular}
\caption{The Randi\'{c} index and harmonic index of real networks.}\label{tabreal}
\end{center}
\end{table}

\section{Proof of main results}\label{proof}
In this section, we provide the detailed proofs of the main results. Recall that  $A_{ij}=1$ if and only if $\{i,j\}$ is an edge. Then the general Randi\'{c} index in (\ref{randicindex})  and the general sum-connectivity index in (\ref{harindexm}) can be written as
\begin{eqnarray*}\label{rdeq1}
\mathcal{R}_{\alpha}&=&\sum_{1\leq i<j\leq n}A_{ij}d_i^{\alpha}d_j^{\alpha},\\
\chi_{\alpha}&=&\sum_{1\leq i<j\leq n}A_{ij}(d_i+d_j)^{\alpha}.
\end{eqnarray*}
Note that the degrees $d_i$ are not independently and identically distributed. Moreover, $\mathcal{R}_{\alpha}$
and $\chi_{\alpha}$ are non-linear functions of $d_i$. These facts make it a non-trivial task to derive the limits of $\mathcal{R}_{\alpha}$
and $\chi_{\alpha}$ for general $\alpha$. The proof strategy is as follows: (a) use the Taylor expansion to expand $\mathcal{R}_{\alpha}$
or $\chi_{\alpha}$ as a sum of leading term and reminder terms; (b) find the order of the leading term and the reminder terms.

\medskip

\noindent
{\bf Proof of Theorem \ref{theorem:1}:} (I) We prove the result of the general Randi\'{c} index first. For convenience, let 
\begin{equation}\label{rrdeq1}
\mathcal{R}_{-\alpha}=\sum_{1\leq i<j\leq n}A_{ij}d_i^{-\alpha}d_j^{-\alpha}.
\end{equation}
We provide the proof in two cases: $\alpha>-1$ and $\alpha\leq-1$. Denote $\mu_i=\mathbb{E}(d_i)=p_nf_i$. 

Let $\alpha>-1$.
Applying the mean value theorem to the mapping $x\rightarrow x^{-\alpha}$, we have
\[\frac{1}{d_i^{\alpha}}=\frac{1}{\mu_i^{\alpha}}-\alpha\frac{d_i-\mu_i}{X_i^{\alpha+1}},\]
where $d_i\leq X_i\leq\mu_i$ or $\mu_i\leq X_i\leq d_i$. Since $A_{ii}=0$ $(i=1,2,\dots,n)$ and the adjacency matrix $A$ is symmetric, by (\ref{rrdeq1}) one has
\begin{eqnarray}\nonumber
\mathcal{R}_{-\alpha}&=&\frac{1}{2}\sum_{1\leq i,j\leq n}\frac{A_{ij}}{d_i^{\alpha}d_j^{\alpha}}\\ \nonumber
&=&\frac{1}{2}\sum_{1\leq i,j\leq n}\frac{A_{ij}}{\mu_i^{\alpha}\mu_j^{\alpha}}-\frac{\alpha}{2}\sum_{1\leq i,j\leq n}\frac{A_{ij}(d_i-\mu_i)}{X_i^{\alpha+1}\mu_j^{\alpha}}-\frac{\alpha}{2}\sum_{1\leq i,j\leq n}\frac{A_{ij}(d_j-\mu_j)}{X_j^{\alpha+1}\mu_i^{\alpha}}\\ \label{rrdeq2}
&&+\frac{\alpha^2}{2}\sum_{1\leq i,j\leq n}\frac{A_{ij}(d_i-\mu_i)(d_j-\mu_j)}{X_i^{\alpha+1}X_j^{\alpha+1}}.
\end{eqnarray}
Next we show the first term in (\ref{rrdeq2}) is the leading term. To this end, we will find the exact order of the first term and show the remaining terms are of smaller order.

Firstly, we show the first term in (\ref{rrdeq2}) is asymptotically equal to its expectation. 
 By the assumption $\min_{1\leq i,j\leq n}\{f_{ij}\}>\epsilon$, it is clear that $np_n\epsilon\leq\mu_i\leq np_n$ for all $i\in[n]$ and $\epsilon n^2\leq \sum_{1\leq i,j\leq n}f_{ij}\leq n^2$. Note that $A_{ij}\ (1\leq i<j\leq n)$ are independent and $\mathbb{E}(A_{ij})=p_nf_{ij}$. Then
\begin{eqnarray*}
\mathbb{E}\left[\sum_{1\leq i<j\leq n}\frac{A_{ij}-p_nf_{ij}}{\mu_i^{\alpha}\mu_j^{\alpha}}\right]^2&=&\sum_{1\leq i<j\leq n}\mathbb{E}\left[\frac{A_{ij}-p_nf_{ij}}{\mu_i^{\alpha}\mu_j^{\alpha}}\right]^2=O\left(\frac{n^2p_n}{(np_n)^{4\alpha}}\right).
\end{eqnarray*}
By the Markov's inequality, it follows that
\[\left|\sum_{1\leq i<j\leq n}\frac{A_{ij}}{\mu_i^{\alpha}\mu_j^{\alpha}}-\sum_{1\leq i<j\leq n}\frac{p_nf_{ij}}{\mu_i^{\alpha}\mu_j^{\alpha}}\right|=\left|\sum_{1\leq i<j\leq n}\frac{A_{ij}-p_nf_{ij}}{\mu_i^{\alpha}\mu_j^{\alpha}}\right|=O_P\left(\frac{\sqrt{n}\sqrt{np_n}}{(np_n)^{2\alpha}}\right).\]
Then we get
\begin{equation}\label{meanc}
    \sum_{1\leq i<j\leq n}\frac{A_{ij}}{\mu_i^{\alpha}\mu_j^{\alpha}}=\sum_{1\leq i<j\leq n}\frac{p_nf_{ij}}{\mu_i^{\alpha}\mu_j^{\alpha}}+O_P\left(\frac{\sqrt{n}\sqrt{np_n}}{(np_n)^{2\alpha}}\right)=\sum_{1\leq i<j\leq n}\frac{p_nf_{ij}}{\mu_i^{\alpha}\mu_j^{\alpha}}\left(1+O_P\left(\frac{1}{\sqrt{n}\sqrt{np_n}}\right)\right).
\end{equation}

Now we find a bound of the second term in (\ref{rrdeq2}). The idea is to find an upper bound of the expectation of its absolute value and then apply the Markov's inequality to get a bound. Note that
\begin{eqnarray}\nonumber
\mathbb{E}\left[\left|\sum_{1\leq i,j\leq n}\frac{A_{ij}(d_i-\mu_i)}{X_i^{\alpha+1}\mu_j^{\alpha}}\right|\right] \label{rdeq3}
&=&\mathbb{E}\left[\left|\sum_{1\leq i\leq n}\left(\sum_{1\leq j\leq n}\frac{A_{ij}}{\mu_j^{\alpha}}\right)\frac{(d_i-\mu_i)}{X_i^{\alpha+1}}\right|\right] \\
&\leq&\mathbb{E}\left[\sum_{1\leq i\leq n}\left(\sum_{1\leq j\leq n}\frac{A_{ij}}{\mu_j^{\alpha}}\right)\frac{|d_i-\mu_i|}{X_i^{\alpha+1}}\right].
\end{eqnarray}
Let $\delta_n=[\log(np_n)]^{-2}$. Recall that $X_i$ is between $d_i$ and $\mu_i$. If $X_i<\delta_n\mu_i$ and $X_i<d_i$, then $X_i<d_i$ and $X_i<\mu_i$. In this case, $X_i$ can not be between $d_i$ and $\mu_i$. Therefore, $X_i<\delta_n\mu_i$ implies $d_i\leq X_i$. Then $I[X_i<\delta_n\mu_i]\leq I[d_i\leq X_i<\delta_n\mu_i]\leq I[X_i<\delta_n\mu_i]$. Note that $np_n\epsilon\leq\mu_i\leq np_n$ for all $i\in[n]$, then we have
\begin{eqnarray}\nonumber
&&\mathbb{E}\left[\left|\sum_{1\leq i,j\leq n}\frac{A_{ij}(d_i-\mu_i)}{X_i^{\alpha+1}\mu_j^{\alpha}}\right|\right] 
\leq O\left(\frac{1}{( np_n)^{\alpha}}\right)\sum_{1\leq i\leq n}\mathbb{E}\left[\frac{d_i|d_i-\mu_i|}{X_i^{\alpha+1}}\right]\\ \nonumber
&=&O\left(\frac{1}{( np_n)^{\alpha}}\right)\sum_{1\leq i\leq n}\mathbb{E}\left[\frac{d_i|d_i-\mu_i|}{X_i^{\alpha+1}}I[\delta_n\mu_i\leq X_i]\right]\\ \nonumber
&&+O\left(\frac{1}{( np_n)^{\alpha}}\right)\sum_{1\leq i\leq n}\mathbb{E}\left[\frac{d_i|d_i-\mu_i|}{X_i^{\alpha+1}}I[\delta_n\mu_i>X_i]\right],\\  \nonumber
&=&O\left(\frac{1}{( np_n)^{\alpha}}\right)\sum_{1\leq i\leq n}\mathbb{E}\left[\frac{d_i|d_i-\mu_i|}{X_i^{\alpha+1}}I[\delta_n\mu_i\leq X_i]\right]\\ \label{rdeq4}
&&+O\left(\frac{1}{( np_n)^{\alpha}}\right)\sum_{1\leq i\leq n}\mathbb{E}\left[\frac{d_i|d_i-\mu_i|}{X_i^{\alpha+1}}I[d_i\leq X_i<\delta_n\mu_i]\right]. 
\end{eqnarray}

Note that $\alpha>-1$. If $\delta_n\mu_i\leq X_i$, then \[\frac{1}{X_i^{\alpha+1}}\leq \frac{1}{(\delta_n\mu_i)^{\alpha+1}}=O\left(\frac{1}{(\delta_nnp_n)^{\alpha+1}}\right).\] 
Hence we have
\begin{eqnarray}\nonumber
&&\frac{1}{( np_n)^{\alpha}}\sum_{1\leq i\leq n}\mathbb{E}\left[\frac{d_i|d_i-\mu_i|}{X_i^{\alpha+1}}I[\delta_n\mu_i\leq X_i]\right]\\ \nonumber
&\leq& O\left(\frac{1}{(\delta_nnp_n)^{\alpha+1}(np_n)^{\alpha}}\right)\sum_{1\leq i\leq n}\mathbb{E}\left[d_i|d_i-\mu_i|I[\delta_n\mu_i\leq X_i]\right]\\ \label{mreveq1}
&\leq& O\left(\frac{1}{(\delta_nnp_n)^{\alpha+1}(np_n)^{\alpha}}\right)\sum_{1\leq i\leq n}\mathbb{E}\left[d_i|d_i-\mu_i|\right].
\end{eqnarray}
 
By definition, the second moment of degree $d_i$ is equal to
\[\mathbb{E}[d_i^2]=\mathbb{E}\left[\sum_{j\neq k}A_{ij}A_{ik}+\sum_{j}A_{ij}\right]=p_n^2\sum_{j\neq k}f_{ij}f_{ik}+p_n\sum_{j}f_{ij},\]
and $Var(d_i)=\sum_{j\neq i}p_nf_{ij}(1-p_nf_{ij})$, then
by the Cauchy-Schwarz inequality, one has
\begin{eqnarray}\nonumber
\sum_{1\leq i\leq n}\mathbb{E}\left[d_i|d_i-\mu_i|\right]
&\leq&\sum_{1\leq i\leq n}\sqrt{\mathbb{E}[d_i^2]\mathbb{E}[(d_i-\mu_i)^2]}\\ \nonumber
&=&\sum_{1\leq i\leq n}\sqrt{\left(p_n^2\sum_{j\neq k}f_{ij}f_{ik}+p_n\sum_{j}f_{ij}\right)\sum_{j}p_nf_{ij}(1-p_nf_{ij})}\\ \label{mreveq2}
&=&O\left(n\sqrt{n^3p_n^3}\right).
\end{eqnarray}
Combining (\ref{mreveq1}) and (\ref{mreveq2}) yields
\begin{eqnarray}\nonumber
\frac{1}{( np_n)^{\alpha}}\sum_{1\leq i\leq n}\mathbb{E}\left[\frac{d_i|d_i-\mu_i|}{X_i^{\alpha+1}}I[\delta_n\mu_i\leq X_i]\right]
&=&O\left(\frac{n\sqrt{n^3p_n^3}}{(\delta_nnp_n)^{\alpha+1}(np_n)^{\alpha}}\right)\\ \nonumber
&=&\frac{n^2p_n}{(np_n)^{2\alpha}}O\left(\frac{1}{\delta_n^{\alpha+1}\sqrt{np_n}}\right)\\ \label{rdeq6}
&=&\frac{n^2p_n}{(np_n)^{2\alpha}}O\left(\frac{(\log (np_n))^{2(\alpha+1)}}{\sqrt{np_n}}\right).
\end{eqnarray}

 Now we bound the second term of (\ref{rdeq4}). 
Note that if $d_i\leq X_i<\delta_n\mu_i$, then $d_i<\mu_i$ and $\frac{d_i}{X_i^{\alpha+1}}\leq \frac{1}{d_i^{\alpha}}$. Since $d_i$ is the degree of node $i$, it can only take integer value between 0 and $n-1$. Moreover, $d_i=0$ implies $A_{ij}=0$ for any $j\in [n]$. By the definition of the Randi\'{c} index
 (\ref{randicindex}), these terms with $d_i=0$ are zero in (\ref{rrdeq1}) and (\ref{rrdeq2}). Therefore, we only consider the terms with $d_i\geq1$ and $d_j\geq1$.
 Then the second term of (\ref{rdeq4}) can be bounded by
\begin{eqnarray}\nonumber
&&\frac{1}{( np_n)^{\alpha}}\sum_{1\leq i\leq n}\mathbb{E}\left[\frac{d_i|d_i-\mu_i|}{X_i^{\alpha+1}}I[d_i\leq X_i<\delta_n\mu_i]\right]\leq\frac{1}{(np_n)^{\alpha}} \sum_{1\leq i\leq n}\mathbb{E}\left[\frac{\mu_i-d_i}{d_i^{\alpha}}I[d_i<\delta_n\mu_i]\right]\\ \label{rheq1}
&=&\frac{1}{( np_n)^{\alpha}}\sum_{1\leq i\leq n}\sum_{k=1}^{\delta_n\mu_i}\frac{\mu_i-k}{k^{\alpha}}\mathbb{P}(d_i=k).
\end{eqnarray}

 Next we obtain an upper bound of $\mathbb{P}(d_i=k)$.
Note that 
the degree $d_i$ follows the Poisson-Binomial distribution $PB(p_nf_{i1},p_nf_{i2},\dots,p_nf_{in})$. Then 
\begin{eqnarray}\nonumber
\mathbb{P}(d_i=k)&=&\sum_{S\subset[n]\setminus\{i\},|S|=k}\prod_{j\in S}p_nf_{ij}\prod_{j\in S^C\setminus\{i\}}(1-p_nf_{ij})\\ \nonumber
&\leq&\sum_{S\subset[n]\setminus\{i\},|S|=k}\prod_{j\in S}p_n\prod_{j\in S^C\setminus\{i\}}(1-p_n\epsilon)\\ \label{rheq2}
&=&\binom{n}{k}p_n^k(1-p_n\epsilon)^{n-k}.
\end{eqnarray}
Note that $\binom{n}{k}\leq e^{k\log n-k\log k+k}$ and 
$(1-p_n\epsilon)^{n-k}=e^{(n-k)\log(1-p_n\epsilon)}$.
Then by (\ref{rheq2}) we get
\begin{eqnarray}\label{rheq3}
\mathbb{P}(d_i=k)&\leq&\exp\left(k\log (np_n)-k\log k+k+(n-k)\log(1-p_n\epsilon)\right).
\end{eqnarray}
Let $g(k)=k\log (np_n)-k\log k+k+(n-k)\log(1-p_n\epsilon)$. Then
\[g^{\prime}(k)=\log\left(\frac{np_n}{1-p_n\epsilon}\right)-\log k.\]
For $k<\frac{np_n}{1-p_n\epsilon}$, $g^{\prime}(k)<0$. For $k>\frac{np_n}{1-p_n\epsilon}$, $g^{\prime}(k)>0$. Hence $g(k)$ achieves its maximum at $k=\frac{np_n}{1-p_n\epsilon}$. For $k\leq \delta_nnp_n$,  $g(k)\leq g(\delta_nnp_n)$. Hence
\[\mathbb{P}(d_i=k)\leq\exp\left(\delta_nnp_n\log\frac{1}{\delta_n(1-p_n\epsilon)}+\delta_nnp_n+n\log(1-p_n\epsilon)\right)\leq \exp\left(-np_n\epsilon(1+o(1))\right).\]

Note that $\mu_i\leq np_n$. Then for $k\leq \delta_n\mu_i\leq \delta_nnp_n$, by (\ref{rheq1}),  (\ref{rheq2}), (\ref{rheq3}), one has
\begin{eqnarray}\nonumber
\mathbb{E}\left[\frac{d_i|d_i-\mu_i|}{X_i^{\alpha+1}}I[d_i\leq X_i<\delta_n\mu_i]\right]
&\leq& \exp\left(\log (\delta_nnp_n)\right)\exp\left(\log (np_n)\right)\exp\left(-np_n\epsilon(1+o(1))\right)\\ \label{mreveq6}
&=&\exp\left(-np_n\epsilon(1+o(1))\right).
\end{eqnarray}
Hence, we get
\begin{eqnarray}\nonumber
\frac{1}{( np_n)^{\alpha}}\sum_{1\leq i\leq n}\mathbb{E}\left[\frac{d_i|d_i-\mu_i|}{X_i^{\alpha+1}}I[d_i\leq X_i<\delta_n\mu_i]\right]
&=&\frac{1}{( np_n)^{\alpha}}ne^{-\epsilon np_n(1+o(1))}=\frac{n^2p_n}{(np_n)^{2\alpha}}e^{-\epsilon np_n(1+o(1))}.\\ \label{mreveq10}
\end{eqnarray}

Recall that $np_n\log 2\geq \log n$. Then $\frac{(\log (np_n))^{s}}{(np_n)^k}e^{-\epsilon np_n(1+o(1))}=o(1)$ for any fixed positive constants $k,s,\epsilon$. By (\ref{rdeq3}), (\ref{rdeq4}), (\ref{rdeq6}), (\ref{mreveq10}) and the Markov's inequality, one has
\begin{eqnarray}\label{nmreveq1}
\sum_{1\leq i,j\leq n}\frac{A_{ij}(d_i-\mu_i)}{X_i^{\alpha+1}\mu_j^{\alpha}}=O_P\left(\frac{n^2p_n}{(np_n)^{2\alpha}}\frac{(\log (np_n))^{2(\alpha+1)}}{\sqrt{np_n}}\right).
\end{eqnarray}

 The third term in (\ref{rrdeq2}) can be similarly bounded as the second term. 
Now we consider the last term in (\ref{rrdeq2}). Note that
\begin{eqnarray}\nonumber
\sum_{1\leq i,j\leq n}\frac{A_{ij}(d_i-\mu_i)(d_j-\mu_j)}{X_i^{\alpha+1}X_j^{\alpha+1}}
&=&
\sum_{1\leq i,j\leq n}\frac{A_{ij}(d_i-\mu_i)(d_j-\mu_j)}{X_i^{\alpha+1}X_j^{\alpha+1}}I[X_i\geq\delta_n\mu_i,X_j\geq\delta_n\mu_j]\\ \nonumber
&&+\sum_{1\leq i,j\leq n}\frac{A_{ij}(d_i-\mu_i)(d_j-\mu_j)}{X_i^{\alpha+1}X_j^{\alpha+1}}I[X_i<\delta_n\mu_i,X_j\geq\delta_n\mu_j]\\ \nonumber
&&+\sum_{1\leq i,j\leq n}\frac{A_{ij}(d_i-\mu_i)(d_j-\mu_j)}{X_i^{\alpha+1}X_j^{\alpha+1}}I[X_i\geq\delta_n\mu_i,X_j<\delta_n\mu_j]\\ \nonumber
&&+\sum_{1\leq i,j\leq n}\frac{A_{ij}(d_i-\mu_i)(d_j-\mu_j)}{X_i^{\alpha+1}X_j^{\alpha+1}}I[X_i<\delta_n\mu_i,X_j<\delta_n\mu_j].\\
\label{mreveq3}
\end{eqnarray}
 We shall bound each term in (\ref{mreveq3}). The first term can be bounded as follows.  
\begin{eqnarray}\nonumber
&&\mathbb{E}\left[\left|\sum_{1\leq i,j\leq n}\frac{A_{ij}(d_i-\mu_i)(d_j-\mu_j)}{X_i^{\alpha+1}X_j^{\alpha+1}}I[X_i\geq\delta_n\mu_i,X_j\geq\delta_n\mu_j]\right|\right]\\ \nonumber
&\leq&\frac{1}{\delta_n^{2(\alpha+1)}}\sum_{1\leq i,j\leq n}\mathbb{E}\left[\frac{A_{ij}|d_i-\mu_i||d_j-\mu_j|}{\mu_i^{\alpha+1}\mu_j^{\alpha+1}}I[X_i\geq\delta_n\mu_i,X_j\geq\delta_n\mu_j]\right]\\  \label{radeq1}
&\leq&\frac{1}{\delta_n^{2(\alpha+1)}}O\left(\frac{1}{(np_n)^{2(\alpha+1)}}\right)\sum_{1\leq i,j\leq n}\mathbb{E}\left[A_{ij}|d_i-\mu_i||d_j-\mu_j|\right]. 
\end{eqnarray}
Denote $\tilde{d}_i=\sum_{k\neq j,i}A_{ik}$, $\tilde{d}_j=\sum_{k\neq j,i}A_{jk}$, $\tilde{\mu}_i=\mathbb{E}(\tilde{d}_i)$ and $\tilde{\mu}_j=\mathbb{E}(\tilde{d}_j)$. Then $\tilde{d}_i$ and $\tilde{d}_j$ are independent, $d_i=\tilde{d}_i+A_{ij}$ and $d_j=\tilde{d}_j+A_{ij}$. It is easy to get that
\[|d_i-\mu_i|=|\tilde{d}_i-\tilde{\mu}_i+A_{ij}-p_nf_{ij}|\leq |\tilde{d}_i-\tilde{\mu}_i|+|A_{ij}-p_nf_{ij}|\leq |\tilde{d}_i-\tilde{\mu}_i|+1,\]
\[\mathbb{E}[|\tilde{d}_i-\tilde{\mu}_i|]\leq \sqrt{\mathbb{E}[(\tilde{d}_i-\tilde{\mu}_i)^2]}= \sqrt{\sum_{k\neq j,i}p_nf_{ik}(1-p_nf_{ik})}=O(\sqrt{np_n}).\]
Similarly, $|d_j-\mu_j|\leq |\tilde{d}_j-\tilde{\mu}_j|+1$ and $\mathbb{E}[|\tilde{d}_j-\tilde{\mu}_j|]=O(\sqrt{np_n})$.
Then we have
\begin{eqnarray}\nonumber
\mathbb{E}\left[A_{ij}|d_i-\mu_i||d_j-\mu_j|\right]&\leq& \mathbb{E}[A_{ij}]+\mathbb{E}[A_{ij}|\tilde{d}_i-\tilde{\mu}_i||\tilde{d}_j-\tilde{\mu}_j|]\\ \nonumber
&&+\mathbb{E}[A_{ij}|\tilde{d}_i-\tilde{\mu}_i|]+\mathbb{E}[A_{ij}|\tilde{d}_j-\tilde{\mu}_j|]\\ \nonumber
&=&p_nf_{ij}+p_nf_{ij}\mathbb{E}[|\tilde{d}_i-\tilde{\mu}_i|]\mathbb{E}[|\tilde{d}_j-\tilde{\mu}_j|]\\  \nonumber
&&+p_nf_{ij}\mathbb{E}[|\tilde{d}_i-\tilde{\mu}_i|]+p_nf_{ij}\mathbb{E}[|\tilde{d}_j-\tilde{\mu}_j|]\\ \label{mreveq4}
&=&O\left(np_n^2\right).
\end{eqnarray}
Combining (\ref{radeq1}) and (\ref{mreveq4}) yields
\begin{eqnarray}\nonumber
&&\mathbb{E}\left[\left|\sum_{1\leq i,j\leq n}\frac{A_{ij}(d_i-\mu_i)(d_j-\mu_j)}{X_i^{\alpha+1}X_j^{\alpha+1}}I[X_i\geq\delta_n\mu_i,X_j\geq\delta_n\mu_j]\right|\right]\\ \nonumber
&\leq&\frac{1}{\delta_n^{2(\alpha+1)}}O\left(\frac{n^3p_n^2}{(np_n)^{2(\alpha+1)}}\right)\\ \nonumber 
&=&\frac{n^2p_n}{(np_n)^{2\alpha}}O\left(\frac{1}{\delta_n^{2(\alpha+1)}np_n}\right)\\ \label{radpr1}
&=&\frac{n^2p_n}{(np_n)^{2\alpha}}O\left(\frac{(\log(np_n))^{4(\alpha+1)}}{np_n}\right),
\end{eqnarray}

The second term in (\ref{mreveq3}) can be bounded as follows.
\begin{eqnarray}\nonumber
&&\mathbb{E}\left[\left|\sum_{1\leq i,j\leq n}\frac{A_{ij}(d_i-\mu_i)(d_j-\mu_j)}{X_i^{\alpha+1}X_j^{\alpha+1}}I[X_i\geq\delta_n\mu_i,d_j\leq X_j<\delta_n\mu_j]\right|\right]\\ \nonumber
&\leq&\frac{1}{\delta_n^{\alpha+1}}\sum_{1\leq i,j\leq n}\mathbb{E}\left[\frac{A_{ij}|d_i-\mu_i||d_j-\mu_j|}{\mu_i^{\alpha+1}d_j^{\alpha+1}}I[X_i\geq\delta_n\mu_i,d_j\leq X_j<\delta_n\mu_j]\right]\\ \label{mreveq6}
&\leq&\frac{1}{\delta_n^{\alpha+1}(np_n)^{\alpha+1}}\sum_{1\leq i,j\leq n}\mathbb{E}\left[\frac{A_{ij}|d_i-\mu_i||d_j-\mu_j|}{d_j^{\alpha+1}}I[d_j<\delta_n\mu_j]\right].
\end{eqnarray}
Recall that 
\[|d_i-\mu_i|=|\tilde{d}_i-\tilde{\mu}_i+A_{ij}-p_nf_{ij}|,\ \ \ \ |d_j-\mu_j|=|\tilde{d}_j-\tilde{\mu}_j+A_{ij}-p_nf_{ij}|.\]
Moreover, $d_j<\delta_n\mu_j$ implies $\tilde{d}_j<\delta_n\mu_j$. Then we have
\begin{eqnarray}\nonumber
&&\mathbb{E}\left[\frac{A_{ij}|d_i-\mu_i||d_j-\mu_j|}{d_j^{\alpha+1}}I[d_j<\delta_n\mu_j]\right]\\ \nonumber
&=&\mathbb{E}\left[\frac{A_{ij}|\tilde{d}_i-\tilde{\mu}_i+A_{ij}-p_nf_{ij}||\tilde{d}_j-\tilde{\mu}_j+A_{ij}-p_nf_{ij}|}{d_j^{\alpha+1}}I[d_j<\delta_n\mu_j]\Big|A_{ij}=1\right] \mathbb{P}(A_{ij}=1)\\ \label{mreveq7}
&\leq&p_n\mathbb{E}\left[\frac{|\tilde{d}_i-\tilde{\mu}_i+1-p_nf_{ij}||\tilde{d}_j-\tilde{\mu}_j+1-p_nf_{ij}|}{(\tilde{d}_j+1)^{\alpha+1}}I[\tilde{d}_j<\delta_n\mu_j]\right].   
\end{eqnarray}
Since $\tilde{d}_i$, $\tilde{d}_j$ are independent and $\mathbb{E}[|\tilde{d}_j-\tilde{\mu}_j|]=O(\sqrt{np_n})$, then by a similar argument as in (\ref{rheq1})-(\ref{mreveq10}), it follows that 
\begin{eqnarray}\nonumber
&&p_n\mathbb{E}\left[\frac{|\tilde{d}_i-\tilde{\mu}_i+1-p_nf_{ij}||\tilde{d}_j-\tilde{\mu}_j+1-p_nf_{ij}|}{(\tilde{d}_j+1)^{\alpha+1}}I[\tilde{d}_j<\delta_n\mu_j]\right]\\ \nonumber
&\leq&p_n\sqrt{np_n}\mathbb{E}\left[\frac{|\tilde{d}_j-\tilde{\mu}_j+1-p_nf_{ij}|}{(\tilde{d}_j+1)^{\alpha+1}}I[\tilde{d}_j<\delta_n\mu_j]\right]\\ \label{mreveq8}
&\leq&p_n\sqrt{np_n}e^{-\epsilon np_n(1+o(1))}.
\end{eqnarray}

Combining (\ref{mreveq6}), (\ref{mreveq7}) and (\ref{mreveq8}) yields

\begin{eqnarray}\nonumber
&&\mathbb{E}\left[\left|\sum_{1\leq i,j\leq n}\frac{A_{ij}(d_i-\mu_i)(d_j-\mu_j)}{X_i^{\alpha+1}X_j^{\alpha+1}}I[X_i\geq\delta_n\mu_i,d_j\leq X_j<\delta_n\mu_j]\right|\right]\\ \nonumber
&\leq&\frac{p_n\sqrt{np_n}}{\delta_n^{\alpha+1}(np_n)^{\alpha+1}}n^2e^{-\epsilon np_n(1+o(1))}\\ \label{nrdeq1}
&=&\frac{n^2p_n}{(np_n)^{2\alpha}}e^{-\epsilon np_n(1+o(1))}.
\end{eqnarray}

The third term in (\ref{mreveq3}) can be similarly bounded as the second term. Now we consider the last term in (\ref{mreveq3}).
By a similar argument as in (\ref{mreveq6})-(\ref{nrdeq1}), one gets

\begin{eqnarray}\nonumber
&&\mathbb{E}\left[\left|\sum_{1\leq i,j\leq n}\frac{A_{ij}(d_i-\mu_i)(d_j-\mu_j)}{X_i^{\alpha+1}X_j^{\alpha+1}}I[d_i\leq X_i<\delta_n\mu_i,d_j\leq X_j<\delta_n\mu_j]\right|\right]\\ \nonumber
&\leq&\sum_{1\leq i,j\leq n}\mathbb{E}\left[\frac{A_{ij}|d_i-\mu_i||d_j-\mu_j|}{d_i^{\alpha+1}d_j^{\alpha+1}}I[d_i\leq\delta_n\mu_i,d_j\leq \delta_n\mu_j]\right]\\ \nonumber
&\leq&\sum_{1\leq i,j\leq n}\mathbb{E}\left[\frac{A_{ij}|\tilde{d}_i-\tilde{\mu}_i+A_{ij}-p_nf_{ij}||\tilde{d}_j-\tilde{\mu}_j+A_{ij}-p_nf_{ij}|}{(\tilde{d}_j+A_{ij})^{\alpha+1}(\tilde{d}_j+A_{ij})^{\alpha+1}}I[\tilde{d}_i\leq\delta_n\mu_i,\tilde{d}_j\leq \delta_n\mu_j]\right]\\ \nonumber
\end{eqnarray}
\begin{eqnarray}\nonumber
&\leq&p_n\sum_{1\leq i,j\leq n}\mathbb{E}\left[\frac{(|\tilde{d}_i-\tilde{\mu}_i|+1)(|\tilde{d}_j-\tilde{\mu}_j|+1)}{(\tilde{d}_j+1)^{\alpha+1}(\tilde{d}_j+1)^{\alpha+1}}I[\tilde{d}_i\leq\delta_n\mu_i,\tilde{d}_j\leq \delta_n\mu_j]\right]\\ \nonumber
&=&p_n\left(\sum_{1\leq i\leq n}\mathbb{E}\left[\frac{(|\tilde{d}_i-\tilde{\mu}_i|+1)}{(\tilde{d}_i+1)^{\alpha+1}}I[\tilde{d}_i\leq\delta_n\mu_i\right]\right)^2\\ \label{nrdeq2}
&\leq&p_nn^2e^{-2\epsilon np_n(1+o(1))}=\frac{n^2p_n}{(np_n)^{2\alpha}}e^{-2\epsilon np_n(1+o(1))}.
\end{eqnarray}

By (\ref{mreveq3})-(\ref{nrdeq2}) and the Markov's inequality, it follows that
\begin{eqnarray}\label{nmreveq2}
\sum_{1\leq i,j\leq n}\frac{A_{ij}(d_i-\mu_i)(d_j-\mu_j)}{X_i^{\alpha+1}X_j^{\alpha+1}}
=O_P\left(\frac{n^2p_n}{(np_n)^{2\alpha}}\frac{(\log(np_n))^{4(\alpha+1)}}{np_n}\right).
\end{eqnarray}

 It is easy to verify that $\sum_{1\leq i<j\leq n}\frac{p_nf_{ij}}{\mu_i^{\alpha}\mu_j^{\alpha}} \geq \frac{\epsilon n(n-1)p_n}{2(np_n)^{2\alpha}}$. Then combining (\ref{rrdeq2}), (\ref{meanc}), (\ref{nmreveq1}) and (\ref{nmreveq2}) yields the limit of $\mathcal{R}_{-\alpha}$ with $\alpha>-1$.

Next, we consider $\mathcal{R}_{-\alpha}$ for $\alpha\leq-1$. In this case, we rewrite the general Randi\'{c} index as
\begin{equation}\label{rdeq1}
\mathcal{R}_{\alpha}=\sum_{1\leq i<j\leq n}A_{ij}d_i^{\alpha}d_j^{\alpha}, \hskip 1cm \alpha\geq1.
\end{equation}
By the Taylor expansion, we have
\[d_i^{\alpha}=\mu_i^{\alpha}+\alpha X_i^{\alpha-1}(d_i-\mu_i),\]
where $X_i$ is between $d_i$ and $\mu_i$.
Then
\begin{eqnarray}\nonumber
\mathcal{R}_{\alpha}&=&\frac{1}{2}\sum_{1\leq i,j\leq n}A_{ij}d_i^{\alpha}d_j^{\alpha}\\ \nonumber
&=&\frac{1}{2}\sum_{1\leq i,j\leq n}A_{ij}\mu_i^{\alpha}\mu_j^{\alpha}+\frac{\alpha}{2}\sum_{1\leq i,j\leq n}A_{ij}(d_i-\mu_i)X_i^{\alpha-1}\mu_j^{\alpha}+\frac{\alpha}{2}\sum_{1\leq i,j\leq n}A_{ij}(d_j-\mu_j)X_j^{\alpha-1}\mu_i^{\alpha}\\ \label{rdeq2}
&&+\frac{\alpha^2}{2}\sum_{1\leq i,j\leq n}A_{ij}(d_i-\mu_i)(d_j-\mu_j)X_i^{\alpha-1}X_j^{\alpha-1}.
\end{eqnarray}

We shall show that the first term in (\ref{rdeq2}) is the leading term and the remaining terms are of smaller order. Similar to (\ref{meanc}), it is easy to get
\begin{equation}\label{meancc}
    \sum_{1\leq i<j\leq n}A_{ij}\mu_i^{\alpha}\mu_j^{\alpha}=\sum_{1\leq i<j\leq n}p_nf_{ij}\mu_i^{\alpha}\mu_j^{\alpha}\left(1+O_P\left(\frac{1}{\sqrt{n}\sqrt{np_n}}\right)\right).
\end{equation}

Since the  second term and the third term in (\ref{rdeq2}) have the same order, we only need to bound the second term and the last term.  Let $M=\frac{4}{\epsilon(1-p_n\epsilon)}$. Clearly $M$ is bounded and $M>4$. The expectation of the absolute value of the second term in (\ref{rdeq2}) can be bounded by
\begin{eqnarray}\nonumber
\mathbb{E}\left[\left|\sum_{1\leq i,j\leq n}A_{ij}(d_i-\mu_i)X_i^{\alpha-1}\mu_j^{\alpha}\right|\right]&\leq&\mathbb{E}\left[\sum_{1\leq i,j\leq n}A_{ij}\left|d_i-\mu_i\right|X_i^{\alpha-1}\mu_j^{\alpha}I[M\mu_i\leq X_i\leq d_i]\right]\\ \label{neweq1}
&&+\mathbb{E}\left[\sum_{1\leq i,j\leq n}A_{ij}\left|d_i-\mu_i\right|X_i^{\alpha-1}\mu_j^{\alpha}I[X_i\leq M\mu_i]\right].
\end{eqnarray}

Note that
\begin{eqnarray}\nonumber
\mathbb{E}\left[\sum_{1\leq i,j\leq n}A_{ij}\left|d_i-\mu_i\right|X_i^{\alpha-1}\mu_j^{\alpha}I[X_i\leq M\mu_i]\right]&\leq& M^{\alpha-1}(np_n)^{2\alpha-1}\sum_{1\leq i,j\leq n}\mathbb{E}\left[A_{ij}\left|\tilde{d}_i-\mu_i+A_{ij}\right|\right]\\
&=&(np_n)^{2\alpha}n^2p_nO\left(\frac{1}{\sqrt{np_n}}\right),
\end{eqnarray}
and
\begin{eqnarray}\nonumber
&&\mathbb{E}\left[\sum_{1\leq i,j\leq n}A_{ij}\left|d_i-\mu_i\right|X_i^{\alpha-1}\mu_j^{\alpha}I[M\mu_i\leq X_i\leq d_i]\right]\\ \nonumber
&\leq&O((np_n)^{\alpha})\mathbb{E}\left[\sum_{1\leq i,j\leq n}A_{ij}\left|d_i-\mu_i\right|d_i^{\alpha-1}I[M\mu_i\leq d_i]\right]\\ \nonumber
&=&O((np_n)^{\alpha}p_n)\sum_{1\leq i,j\leq n}\mathbb{E}\left[\left|\tilde{d}_i-\tilde{\mu}_i+1-p_nf_{ij}\right|\tilde{d}_i^{\alpha-1}I[M\mu_i-1\leq \tilde{d}_i]\right]\\ \label{mreveq13}
&=&O((np_n)^{\alpha}p_n)\sum_{1\leq i,j\leq n}\sum_{k=M\mu_i-1}^{n-2}k^{\alpha-1}(k-\tilde{\mu}_i+1-p_nf_{ij})\mathbb{P}(\tilde{d}_i=k).
\end{eqnarray}
By a similar argument as in (\ref{rheq3}),  it follows that
\begin{eqnarray}\nonumber
\sum_{k=M\mu_i-1}^{n-2}k^{\alpha-1}(k-\tilde{\mu}_i+1-p_nf_{ij})\mathbb{P}(\tilde{d}_i=k)
&\leq& \sum_{k=M\mu_i-1}^{n-2}k^{\alpha}\binom{n}{k}p_n^k(1-p_n\epsilon)^{n-k}\\
&\leq&\sum_{k=M\mu_i-1}^{n-2}\exp\left(\alpha\log k+g(k)\right).
\end{eqnarray}
Let $h(k)=\alpha\log k+g(k)$. Then
\[h^{\prime}(k)=\frac{\alpha}{k}+\log\left(\frac{np_n}{1-p_n\epsilon}\right)-\log k.\]
Hence $h(k)$ is decreasing for $k>\frac{1.1np_n}{1-p_n\epsilon}$ and large $n$. Since $k\geq M\mu_i-1\geq M\epsilon np_n-1\geq\frac{2np_n}{1-p_n\epsilon}$ for large $n$, then 
\[h(k)\leq h\left(\frac{2np_n}{1-p_n\epsilon}\right)=\alpha\log\left(\frac{2np_n}{1-p_n\epsilon}\right)-\frac{2np_n\log 2}{1-p_n\epsilon}+n\log(1-p_n\epsilon)\leq -\frac{np_n\log 2}{1-p_n\epsilon}-\epsilon np_n.\]

By the assumption  $np_n\log 2\geq \log n$, it is easy to get $\log n-\frac{np_n\log 2}{1-p_n\epsilon}<0$. Then
\begin{eqnarray}\nonumber
\sum_{k=M\mu_i-1}^{n-2}k^{\alpha-1}(k-\tilde{\mu}_i+1-p_nf_{ij})\mathbb{P}(\tilde{d}_i=k)
&\leq&n\exp\left(-\frac{np_n\log 2}{1-p_n\epsilon}-\epsilon np_n\right)\\ \label{expeq1}
&\leq&\exp\left(-\epsilon np_n(1+o(1))\right).
\end{eqnarray}
Hence (\ref{neweq1}) is bounded by  $(np_n)^{2\alpha}n^2p_nO\left(\frac{1}{\sqrt{np_n}}\right)$.

Now we bound the last term in (\ref{rdeq2}). Note that
\begin{eqnarray}\nonumber
&&\sum_{1\leq i,j\leq n}|A_{ij}(d_i-\mu_i)(d_j-\mu_j)X_i^{\alpha-1}X_j^{\alpha-1}|\\ \nonumber
&=&\sum_{1\leq i,j\leq n}|A_{ij}(d_i-\mu_i)(d_j-\mu_j)X_i^{\alpha-1}X_j^{\alpha-1}|I[X_i\leq M\mu_i,X_j\leq M\mu_j]\\ \nonumber
&&+\sum_{1\leq i,j\leq n}|A_{ij}(d_i-\mu_i)(d_j-\mu_j)X_i^{\alpha-1}X_j^{\alpha-1}|I[X_i\leq M\mu_i,X_j\geq M\mu_j]\\ \nonumber
&&+\sum_{1\leq i,j\leq n}|A_{ij}(d_i-\mu_i)(d_j-\mu_j)X_i^{\alpha-1}X_j^{\alpha-1}|I[X_i\geq M\mu_i,X_j\leq M\mu_j]\\ \label{mreveq11}
&&+\sum_{1\leq i,j\leq n}|A_{ij}(d_i-\mu_i)(d_j-\mu_j)X_i^{\alpha-1}X_j^{\alpha-1}|I[X_i\geq M\mu_i,X_j\geq M\mu_j].
\end{eqnarray}
Since $X_i$ is between $d_i$ and $\mu_i$, then $X_i\leq M\mu_i$ implies $d_i\leq X_i\leq M\mu_i$, and $X_i\geq M\mu_i$ implies $d_i\geq X_i\geq M\mu_i$. Similar results hold for $X_j$. Then by (\ref{mreveq11}) we have 
\begin{eqnarray}\nonumber
&&\sum_{1\leq i,j\leq n}|A_{ij}(d_i-\mu_i)(d_j-\mu_j)X_i^{\alpha-1}X_j^{\alpha-1}|\\ \nonumber
&\leq&\sum_{1\leq i,j\leq n}|A_{ij}(d_i-\mu_i)(d_j-\mu_j)X_i^{\alpha-1}X_j^{\alpha-1}|I[X_i\leq M\mu_i, X_j\leq M\mu_j]\\ \nonumber
&&+\sum_{1\leq i,j\leq n}|A_{ij}(d_i-\mu_i)(d_j-\mu_j)X_i^{\alpha-1}X_j^{\alpha-1}|I[X_i\leq M\mu_i,d_j\geq X_j\geq M\mu_j]\\ \nonumber
&&+\sum_{1\leq i,j\leq n}|A_{ij}(d_i-\mu_i)(d_j-\mu_j)X_i^{\alpha-1}X_j^{\alpha-1}|I[d_i\geq X_i\geq M\mu_i,X_j\leq M\mu_j]\\ \label{mreveq12}
&&+\sum_{1\leq i,j\leq n}|A_{ij}(d_i-\mu_i)(d_j-\mu_j)X_i^{\alpha-1}X_j^{\alpha-1}|I[d_i\geq X_i\geq M\mu_i,d_j\geq X_j\geq M\mu_j].
\end{eqnarray}
Now we bound the expectation of each term in (\ref{mreveq12}). Since the second term and the third term have the same order, it suffices to bound the first term, second term and the last term. By a similar argument as in (\ref{mreveq13}) and (\ref{expeq1}), it is easy to get the following results.
\begin{eqnarray}\nonumber
&&\mathbb{E}\left[\sum_{1\leq i,j\leq n}A_{ij}|d_i-\mu_i||d_j-\mu_j|X_i^{\alpha-1}X_j^{\alpha-1}I[X_i\leq M\mu_i,X_j\leq M\mu_j]\right]\\ \nonumber
&\leq&O((np_n)^{2(\alpha-1)}p_n)\sum_{1\leq i,j\leq n}\mathbb{E}|\tilde{d}_i-\tilde{\mu}_i+1-p_nf_{ij}||\tilde{d}_j-\tilde{\mu}_j+1-p_nf_{ij}|\\ \nonumber
&=&O((np_n)^{2(\alpha-1)}p_nn^2np_n)\\ \label{mreveq14}
&=&(np_n)^{2\alpha}n^2p_nO\left(\frac{1}{np_n}\right),
\end{eqnarray}

\begin{eqnarray}\nonumber
&&\mathbb{E}\left[\sum_{1\leq i,j\leq n}A_{ij}|d_i-\mu_i||d_j-\mu_j|X_i^{\alpha-1}X_j^{\alpha-1}I[d_i\geq X_i\geq M\mu_i,d_j\geq X_j\geq M\mu_j]\right]\\ \nonumber
&\leq&\mathbb{E}\left[\sum_{1\leq i,j\leq n}A_{ij}d_i^{\alpha}d_j^{\alpha}I[d_i\geq  M\mu_i,d_j\geq M\mu_j]\right]\\ \nonumber
&\leq&p_n\sum_{1\leq i,j\leq n}\mathbb{E}\left[(\tilde{d}_i+1)^{\alpha}(\tilde{d}_j+1)^{\alpha}I[\tilde{d}_i\geq  M\mu_i-1,\tilde{d}_j\geq M\mu_j-1]\right]\\ \nonumber
&=&p_n\left(\sum_{1\leq i\leq n}\mathbb{E}\left[(\tilde{d}_i+1)^{\alpha}I[\tilde{d}_i\geq  M\mu_i-1]\right]\right)^2\\ \label{mreveq15}
&=& O\left(n^2p_n\right)\exp\left(-2\epsilon np_n(1+o(1))\right),
\end{eqnarray}
and
\begin{eqnarray}\nonumber
&&\mathbb{E}\left[\sum_{1\leq i,j\leq n}A_{ij}|d_i-\mu_i||d_j-\mu_j|X_i^{\alpha-1}X_j^{\alpha-1}I[X_i\leq M\mu_i,d_j\geq X_j\geq M\mu_j]\right]\\ \nonumber
&\leq&O((np_n)^{\alpha-1})\mathbb{E}\left[\sum_{1\leq i,j\leq n}A_{ij}|d_i-\mu_i|d_j^{\alpha}I[d_j\geq M\mu_j]\right]\\ \nonumber
&\leq&O((np_n)^{\alpha-1}p_n)\sum_{1\leq i,j\leq n}\mathbb{E}\left[|\tilde{d}_i-\tilde{\mu}_i+1-p_nf_{ij}|(\tilde{d}_j+1)^{\alpha}I[\tilde{d}_j\geq M\mu_j-1]\right]\\ \label{nmreveq16}
&=&O((np_n)^{\alpha-1}p_nn^2\sqrt{np_n})\exp\left(-\epsilon np_n(1+o(1))\right).
\end{eqnarray}
 Combining (\ref{rdeq2})- (\ref{nmreveq16}) yields the desired result. 
Then the proof of the result of the general Randi\'{c} index is complete.

\medskip

(II). Now we prove the result of the general sum-connectivity index.
We provide the proof in two cases: $\alpha<1$ and $\alpha\geq 1$.

Firstly we work on $\chi_{-\alpha}$ with $\alpha>-1$. By Taylor expansion or the mean value theorem, we have
\begin{equation}\label{mreveq16}
\chi_{-\alpha}=\frac{1}{2}\sum_{1\leq i,j\leq n}\frac{A_{ij}}{(d_i+d_j)^{\alpha}}=\frac{1}{2}\sum_{1\leq i,j\leq n}\frac{A_{ij}}{(\mu_i+\mu_j)^{\alpha}}-\frac{\alpha }{2}\sum_{1\leq i,j\leq n}\frac{A_{ij}}{X_{ij}^{\alpha+1}}(d_i-\mu_i+d_j-\mu_j),
\end{equation}
where $X_{ij}$ is between $\mu_i+\mu_j$ and $d_i+d_j$.  We shall prove the first term is the leading term and the second term has smaller order than the first term.

By a similar argument as in (\ref{meanc}), it is easy to get
\begin{equation}\label{harineq3}
\sum_{i<j}\frac{A_{ij}}{(\mu_i+\mu_j)^{\alpha}}=\sum_{i<j}\frac{p_nf_{ij}}{(\mu_i+\mu_j)^{\alpha}}\left(1+O_P\left(\frac{1}{\sqrt{n^2p_n}}\right)\right).
\end{equation}
Hence the first term of (\ref{mreveq16}) is asymptotically equal to $\sum_{i<j}\frac{p_nf_{ij}}{(\mu_i+\mu_j)^{\alpha}}$.

Let $\delta_n=[\log(np_n)]^{-2}$. Since $X_{ij}$ is between $\mu_i+\mu_j$ and $d_i+d_j$, $X_{ij}\leq \delta_n(\mu_i+\mu_j)$ implies $d_i+d_j\leq X_{ij}\leq \delta_n(\mu_i+\mu_j)$. Then
\begin{eqnarray}\nonumber
&&\sum_{i,j}\left|\frac{A_{ij}}{X_{ij}^{\alpha+1}}(d_i-\mu_i+d_j-\mu_j)\right|\\ \nonumber
&\leq&\sum_{i,j}\left|\frac{A_{ij}}{X_{ij}^{\alpha+1}}(d_i-\mu_i+d_j-\mu_j)\right|I[d_i+d_j\leq X_{ij}\leq \delta_n(\mu_i+\mu_j)]\\ \label{mreveq15}
&&+\sum_{i,j}\left|\frac{A_{ij}}{X_{ij}^{\alpha+1}}(d_i-\mu_i+d_j-\mu_j)\right|I[X_{ij}\geq \delta_n(\mu_i+\mu_j)].
\end{eqnarray}
 Next we bound the expectation of each term in (\ref{mreveq15}). For the second term, the expectation can be bounded as follows.
\begin{eqnarray}\nonumber
&&\mathbb{E}\left[\sum_{i,j}\frac{A_{ij}}{X_{ij}^{\alpha+1}}(|d_i-\mu_i|+|d_j-\mu_j|)I[X_{ij}\geq \delta_n(\mu_i+\mu_j)]\right]\\ \nonumber
&\leq&O\left(\frac{1}{\delta_n^{\alpha+1}(np_n)^{\alpha+1}}\right)\sum_{i,j}\mathbb{E}\left[A_{ij}(|\tilde{d}_i-\mu_i+A_{ij}|+|\tilde{d}_j-\mu_i+A_{ij}|)\right]\\ \label{harineq2}
&=&O\left(\frac{n^2p_n\sqrt{np_n}}{\delta_n^{\alpha+1}(np_n)^{\alpha+1}}\right)=\frac{n^2p_n}{(np_n)^{\alpha}}O\left(\frac{[\log (np_n)]^{2(\alpha+1)}}{\sqrt{np_n}}\right).
\end{eqnarray}
Next we focus on the first term in (\ref{mreveq15}). It is clear that
\begin{eqnarray}\nonumber
&&\mathbb{E}\left[\sum_{i,j}\frac{A_{ij}}{X_{ij}^{\alpha+1}}(|d_i-\mu_i|+|d_j-\mu_j|)I[d_i+d_j\leq X_{ij}< \delta_n(\mu_i+\mu_j)]\right]\\ \nonumber
&\leq& \mathbb{E}\left[\sum_{i,j}\frac{A_{ij}(|d_i-\mu_i|+|d_j-\mu_j|)}{(d_i+d_j)^{\alpha+1}}I[d_i+d_j< \delta_n(\mu_i+\mu_j)]\right].
\end{eqnarray}
Note that $d_i+d_j< \delta_n(\mu_i+\mu_j)$ implies $d_i< \delta_n(\mu_i+\mu_j)$ and $d_j< \delta_n(\mu_i+\mu_j)$, and

\[\frac{|d_i-\mu_i|+|d_j-\mu_j|}{(d_i+d_j)^{\alpha+1}}=\frac{|d_i-\mu_i|}{(d_i+d_j)^{\alpha+1}}+\frac{|d_j-\mu_j|}{(d_i+d_j)^{\alpha+1}}\leq\frac{|d_i-\mu_i|}{d_i^{\alpha+1}}+\frac{|d_j-\mu_j|}{d_j^{\alpha+1}}.\]
Then we have
\begin{eqnarray}\nonumber
&&\mathbb{E}\left[\sum_{i,j}\frac{A_{ij}}{X_{ij}^{\alpha+1}}(|d_i-\mu_i|+|d_j-\mu_j|)I[d_i+d_j\leq X_{ij}< \delta_n(\mu_i+\mu_j)]\right]\\ \nonumber
&\leq & \mathbb{E}\left[\sum_{i,j}\frac{A_{ij}|d_i-\mu_i|}{d_i^{\alpha+1}}I[d_i< \delta_n(\mu_i+\mu_j)]\right]+\mathbb{E}\left[\sum_{i,j}\frac{A_{ij}|d_j-\mu_j|}{d_j^{\alpha+1}}I[d_j< \delta_n(\mu_i+\mu_j)]\right]\\ \nonumber
&\leq&2p_n \mathbb{E}\left[\sum_{i,j}\frac{|\tilde{d}_i-\mu_i+1|}{(\tilde{d}_i+1)^{\alpha+1}}I[\tilde{d}_i< \delta_n(\mu_i+\mu_j)]\right]\\ \label{harineq1}
&=&n^2p_ne^{-\epsilon np_n(1+o(1))}=\frac{n^2p_n}{(np_n)^{\alpha}}e^{-\epsilon np_n(1+o(1))}.
\end{eqnarray}
Combining (\ref{mreveq16})- (\ref{harineq1}) yields
\[\chi_{-\alpha}=p_n^{1-\alpha}\sum_{i<j}\frac{f_{ij}}{(f_i+f_j)^{\alpha}}\left(1+O_P\left(\frac{[\log (np_n)]^{2(\alpha+1)}}{\sqrt{n^2p_n}}\right)\right),\hskip 1cm \alpha>-1.\]

Now we work on $\chi_{\alpha}$ with $\alpha\geq 1$. When $\alpha=1$, the proof is trivial. We will focus on $\alpha>1$. By the mean value theorem, one has
\begin{equation}\label{mreveq20}
\chi_{\alpha}=\frac{1}{2}\sum_{i,j}A_{ij}(d_i+d_j)^{\alpha}=\frac{1}{2}\sum_{i,j}A_{ij}(\mu_i+\mu_j)^{\alpha}+\frac{\alpha}{2}\sum_{i,j}A_{ij}X_{ij}^{\alpha-1}(d_i-\mu_i+d_j-\mu_j),
\end{equation}
where $X_{ij}$ is between $\mu_i+\mu_j$ and $d_i+d_j$.

The remaining proof is similar to the proof of the case $\alpha<1$. 
Let $M=\frac{4}{\epsilon(1-p_n\epsilon)}$. It is clear $M$ is bounded and $M>4$.
Note that
\begin{eqnarray}\label{abceq2}
\sum_{i,j}\mathbb{E}\left[A_{ij}X_{ij}^{\alpha-1}(|d_i-\mu_i|+|d_j-\mu_j|)I[X_{ij}\leq M(\mu_i+\mu_j)]\right]=(np_n)^{\alpha}n^2p_nO\left(\frac{1}{\sqrt{np_n}}\right),
\end{eqnarray}
and
\begin{eqnarray}\nonumber
&&\sum_{i,j}\mathbb{E}\left[A_{ij}X_{ij}^{\alpha-1}(|d_i-\mu_i+d_j-\mu_j|)I[d_i+d_j\geq X_{ij}> M(\mu_i+\mu_j)]\right]\\ \nonumber
&\leq&O(1)\sum_{i,j}\mathbb{E}\left[A_{ij}(\tilde{d}_i+\tilde{d}_j+2A_{ij})^{\alpha-1}(|\tilde{d}_i+\tilde{d}_j-\mu_i-\mu_j+2A_{ij}|)I[\tilde{d}_i+\tilde{d}_j>M(\mu_i+\mu_j-1)]\right]\\ \nonumber
&\leq&O(1)p_n\sum_{i,j}\mathbb{E}\left[(\tilde{d}_i+\tilde{d}_j+2)^{\alpha-1}(\tilde{d}_i+\tilde{d}_j)I[\tilde{d}_i+\tilde{d}_j>M(\mu_i+\mu_j-1)]\right]\\ \nonumber
&\leq&O(1)p_n\sum_{i,j}\sum_{k=M(\mu_i+\mu_j-1)}^{2(n-2)}(k+2)^{\alpha-1}k\mathbb{P}(\tilde{d}_i+\tilde{d}_j=k)\\ \label{abceq1}
&=&n^2p_nne^{-\epsilon np_n(1+o(1))}=(np_n)^{\alpha}n^2p_ne^{-\epsilon np_n(1+o(1))},
\end{eqnarray}
where the second last step follows from a similar argument as in (\ref{expeq1}) by noting that  $\tilde{d}_i+\tilde{d}_j$ follows the Poisson-Binomial distribution.

Combining (\ref{mreveq20}), (\ref{abceq2}) and (\ref{abceq1}) yields
\[\chi_{\alpha}=\left(1+O_P\left(\frac{1}{\sqrt{np_n}}\right)\right)p_n^{\alpha+1}\sum_{i<j}(f_i+f_j)^{\alpha}f_{ij},\hskip 1cm \alpha\geq1.\]
Then the proof is complete.

\qed

\section*{Conflict of interest}

The author has no conflict of interest to disclose.

\section*{Acknowledgement}
The author is grateful to the anonymous referees for valuable comments that significantly improve this manuscript.


\begin{thebibliography}{9}
\bibitem{A18}
Abbe, E., Community Detection and Stochastic Block Models:
Recent Developments. {\em Journal of Machine Learning Research}. 2018, {\em 18}, 1-86.



\bibitem{BM05}
Bianconi, G. and Marsili, M. (2005),
Emergence of large cliques in random scale-free network,\textit{Europhysics Letters},74,740.



\bibitem{BM06}
Bianconi, G. and Marsili, M. (2006).
Number of cliques in random scale-free network ensembles,\textit{Physica D: Nonlinear Phenomena}, 224,:1-6.



\bibitem{BCH20}
Bogerd,K., Castro, R., and  Hofstad, R.(2020).
Cliques in rank-1 random graphs: The role of inhomogeneity,\textit{Bernoulli}, 26(1): 253-285 .

\bibitem{BS16}
Bickel, P. J. and Sarkar, P. (2016). Hypothesis testing for automated community detection in networks. 
\textit{Journal of Royal Statistical Society, Series B}, 78, 253-273.


\bibitem{BT78}
Bonchev, D. and Trinajstic, N.(1978). On topological characterization of molecular branching. \textit{International Journal of Quantum Chemistry: Quantum Chemistry Symposium}, 12,293-303.




\bibitem{BDM06}
Britton,T., Deijfen,M.,  Martin-Lof,A.(2006). Generating simple random graphs with prescribed
degree distribution, \textit{Journal of Statistical Physics}, 124:1377–1397.


\bibitem{BE98}
Bollob\'{a}s, B.,  Erdos, P. (1998). Graphs of extremal weights, \textit{Ars Comb.} 50: 225-233.




\bibitem{BES99}
Bollobás B., Erdős P., Sarkar A. (1999).
Extremal graphs for weights
\textit{Discrete Math.}, 200:5-19






\bibitem{CGL16}
 Chiasserini, C.F.,  Garetto, M. and  Leonardi, E. (2016). Social Network De-Anonymization Under
Scale-Free User Relations. \textit{IEEE/ACM Transactions on Networking} 24 (6):3756–3769.



\bibitem{CSN09}
Clauset, A., Shalizi, C. R. and Newman, M.(2009). Power-law Distributions in Empirical Data. \textit{SIAM
review} 51(4), 661–703.




\bibitem{CFK10}
Cavers,M., Fallat, S., Kirkland,S. (2010). On the normalized Laplacian energy and general
Randi\'{c} index $r_1$ of graphs, \textit{Lin. Algebra Appl.}, 433: 172-190.


\bibitem{CHHS21}
Chakrabarty, A., Hazra, S. R., Hollander, F. D. and Sfragara, M.(2021).
Spectra of adjacency and Laplacian matrices of inhomogeneous   Erd\H{o}s-R\'{e}nyi random graphs,
\textit{Random matrices: Theory and applications}, 10(1),215009.

\bibitem{CHHS20}
Chakrabarty, A., Hazra, S. R., Hollander, F. D. and Sfragara, M.(2020).
Large deviation principle for the maximal eigenvalue of
inhomogeneous Erd\H{o}s-R\'{e}nyi random graphs,
\textit{Journal of Theoretical Probability}, https://doi.org/10.1007/s10959-021-01138-w


\bibitem{CCH20}
Chakrabarty, A., Chakrabarty, S. and Hazra, R. S.(2020).
Eigenvalues outside the bulk of  of
inhomogeneous Erd\H{o}s-R\'{e}nyi random graphs,
\textit{Journal of Statistical Physics}, 181: 1746-1780.

\bibitem{CGL16}
 Chiasserini, C.F.,  Garetto, M. and  Leonardi, E. (2016). Social Network De-Anonymization Under
Scale-Free User Relations. \textit{IEEE/ACM Transactions on Networking} 24 (6):3756–3769.

\bibitem{DSG17}
Das K.C., Sun S., Gutman I.(2017).
Normalized Laplacian eigenvalues and Randić energy of graphs
\textit{MATCH Commun. Math. Comput. Chem.}, 77:45-59


\bibitem{DMRSV18}
De Meo, P. etc. (2018).
Estimating Graph Robustness
Through the Randic Index,
\textit{IEEE TRANSACTIONS ON CYBERNETICS}, 48(11):3232-3242.


\bibitem{DMMRSF21}
Dattola, S., etc.(2021). Testing graph robustness indexes for EEG analysis in
Alzheimer’s disease diagnosis,\textit{Electronics},10,1440.



\bibitem{DHHIR20}
Doslic, T., etc.(2020).
On generalized Zagreb indices
of random graphs.
\textit{MATCH Commun. Math. Comput. Chem.} 84:499-511.


\bibitem{E10}
Estrada, E.(2010).
Quantifying network heterogeneity,\textit{PHYSICAL REVIEW E}, 82, 066102.


\bibitem{F87}
 Fajtlowicz, S.(1987).
On conjectures of Graffiti—II
\textit{Congr. Numer.}, 60: 187-197


\bibitem{FMS93}
Favaron, O., Mahéo M. and  Saclé, J.F.(1993).
Some eigenvalue properties in graphs (conjectures of Graffiti—II)
\textit{Discrete Math.}, 111: 197-220

\bibitem{FT13}
Fourches, D. and Tropsha, A.(2013). Using graph indices for the analysis and comparison of chemical datasets. \textit{Molecular Informatics}, 32(9-10), 827-842.


\bibitem{GZFA10}
Goldenberg, A.,etc.(2010). A
survey of statistical network models. \textit{Foundations and Trends
in Machine Learning},
2(2):129–233.


\bibitem{JLN10}
S. Janson, T. Luczak, I. Norros,(2010). Large cliques in a power-law random graph, \textit{J. Appl. Prob.}
47:  1124–1135


\bibitem{JLS19}
A. Janssen, J. Leeuwaarden, S. Shneer, (2019). Counting cliques and cycles in scale-free inhomoge-
neous random graphs, \textit{Journal of Statistical Physics}, 175:161–184.


\bibitem{K07}
Kardar, M.(2007). Statistical physics of particles, Cambridge University Press, Cambridge, 2007.




\bibitem{K09}
Kolaczyk, E.(2009). Statistical analysis of network data. Springer.

\bibitem{LS08}
Li, X. and  Shi, Y.(2008). A survey on the Randic index.
\textit{MATCH Commun. Math. Comput. Chem.}, 59,127-156.


\bibitem{LSG21}
Li,S., Shi, L., and Gao, W.(2021).
Two modified Zagreb indices for random
structures,\textit{Main Group Met. Chem. }. 44: 150–156.

\bibitem{MCSGDF18}
Ma,Y. D., etc.(2018).
From the connectivity index to
various randic-type descriptors.\textit{MATCH Commun. Math. Comput. Chem.}. 80:85-106.


\bibitem{MMRS20}
Martinez-Martinez C.T., Mendez-Bermudez, J.A., Rodriguez, J. and Sigarreta, J.
 (2020).
Computational and analytical studies of the Randi\'{c} index in Erd\H{o}s-R\'{e}nyi models.
\textit{Applied Mathematics and Computation},377,125137.


\bibitem{MMRS21}
Martinez-Martinez C.T., Mendez-Bermudez, J.A., Rodriguez, J. and Sigarreta, J.
 (2021).
Computational and Analytical Studies of
the Harmonic Index on Erd\H{o}s-R\'{e}nyi Models,
\textit{MATCH Commun. Math. Comput. Chem.} 85:395-426.




\bibitem{N03}
Newman, M. (2003). The Structure and Function of Complex Networks. \textit{SIAM review} 45, (2), 167–256


\bibitem{NKMT03}
Nikolic, S, etc.(2003).
The zagreb indices 30 years after. 
\textit{CROATICA CHEMICA ACTA}, 76: 113-124.





\bibitem{N09}
Newman, M.(2009). Networks: an introduction. Oxford University Press, 2009.




\bibitem{R75}
 Randi\'{c} ,M.(1975). Characterization of molecular branching, \textit{J. Am. Chem. Soc.} 97 (23):6609–6615.

\bibitem{RNP16}
Randic M., NoviCM., Plavsic D.(2016).
Solved and Unsolved Problems in Structural Chemistry, CRC Press, Boca Raton.


\bibitem{R08}
Randić M.
On history of the Randić index and emerging hostility toward chemical graph theory.
\textit{MATCH Commun. Math. Comput. Chem.}, 59:5-124.




\bibitem{RS17}
RodrIguez, J. M. and Sigarreta  J. M.(2017). New results on the harmonic index and its
generalizations, 
\textit{MATCH Commun. Math. Comput. Chem.} 78:387-404.



\bibitem{YXL21}
Yu,L., Xu, J. and Lin,X.(2021).
The power of D-hops in matching power-law graphs.\textit{Proceedings of the ACM on Measurement and Analysis of Computing Systems},5(2):1–43.




\bibitem{ZT09}
Zhou,B., Trinajstic, N.(2009). On a novel connectivity index, \textit{J. Math. Chem.}. 46:
1252-1270.

\bibitem{ZT10}
Zhou,B., Trinajstic, N.(2010). On general sum-connectivity index, \textit{J. Math. Chem.} 47:
210-218.

\bibitem{Z12}
Zhong,L.(2012). The harmonic index for graphs, \textit{Appl. Math. Lett.} 25: 561-566.





\end{thebibliography}
\end{document}